\documentclass[12pt]{article}
\usepackage{amssymb}
\usepackage{amscd}
\usepackage{amsmath}
\usepackage[utf8]{inputenc}
\usepackage{amsfonts}
\usepackage{amstext}
\usepackage{color}
\usepackage{a4wide}
\usepackage{graphicx}

\newtheorem{theorem}{Theorem}
\newtheorem{ex}{Example}
\newtheorem{lemma}{Lemma}

\newtheorem{remark}{Remark}
\newtheorem{proposition}{Proposition}
\newtheorem{Def}{Definition}
\newtheorem{corollary}{Corollary}

\newcommand{\La}{\Lambda}
\newcommand{\B}{\mathcal{B}} 
\newcommand{\R}{{\mathbb R}}  \newcommand{\Z}{{\mathbb Z}} \newcommand{\N}{{\mathbb N}}
\newcommand{\T}{{\mathbb T}}
\newcommand{\C}{{\mathbb C}}

\newcommand{\re}{{\rm Re \,}}

\newcommand{\dist}{{\rm dist \,}}
\newcommand{\supp}{{\rm supp \,}}

\newcommand{\ti}{\tilde{t}}

\providecommand{\keywords}[1]
{
  \small	
  \textbf{\textit{Keywords: }} #1
}

\begin{document}
\author{Ilya Zlotnikov}
\title{ On planar sampling with Gaussian kernel in spaces of bandlimited functions\thanks{This research was supported by the Russian Science Foundation (grant No. 18-11-00053), https://rscf.ru/project/18-11-00053/ }}
\maketitle
\begin{abstract}

Let $I=(a,b)\times(c,d)\subset \R_{+}^2$ be an index set and let $\{G_{\alpha}(x) \}_{\alpha \in I}$ be a collection of Gaussian functions,
 i.e. $G_{\alpha}(x) = \exp(-\alpha_1 x_1^2 - \alpha_2 x_2^2)$, where $\alpha = (\alpha_1, \alpha_2) \in I, \, x = (x_1, x_2) \in \R^2$. We present a complete description of the uniformly discrete sets $\Lambda \subset \R^2$ such that every bandlimited signal $f$ admits a stable reconstruction from the samples  $\{f \ast G_{\alpha} (\lambda)\}_{\lambda \in \La}$.

\end{abstract}


\section{Introduction}

The sampling problem deals with recovery of band-limited signals $f$ from the collection of measurements $\{f(\lambda)\}_{\lambda \in \Lambda}$ taken at the points of some uniformly discrete set $\Lambda \subset \R^d$. The classical results deal with one dimensional signals that are elements of the Paley-Wiener or  Bernstein spaces over a fixed interval $[-\sigma,\sigma]$. The sets $\Lambda$ that provide the stable reconstruction, in this case, are completely described. For the Bernstein spaces, the answer is given in terms of a certain density of $\Lambda$ and bandwidth parameter $\sigma$, see \cite{B}. The result for Paley-Wiener spaces is more complicated, see \cite{OS} and \cite{Seip}. It cannot be expressed in terms of a density of $\Lambda$.  We refer the reader to \cite{B} and  \cite{Seip} for the detailed exposition and the proofs.

The complexity of the task significantly increases in the multi-dimensional setting.  Landau \cite{L} proved that the necessary conditions for stable sampling remain valid for the Paley-Wiener spaces over any domain (see \cite{NO} for a much simpler proof). A sufficient condition for a sampling of signals from the Bernstein space with spectrum in a ball was obtained by Beurling, see \cite{B-n}. We also refer the reader to \cite{OUM} for some extensions. However, there is a gap between the necessary and sufficient conditions. Moreover, even for the simplest spectra as balls or cubes, examples show that no description of sampling sets is possible in terms of density of $\Lambda$, see Section 5.7 in~\cite{ou1}.

Recently the so-called dynamical sampling problem (in what follows, we will more often use the term space-time sampling problem) attracted a lot of attention, see \cite{A1}, \cite{A2}, \cite{g}, \cite{UZ}, and references therein.
The dynamical sampling problem  deals with the reconstruction of the initial signal from the given space-time samples.

In this paper, we consider one of the problems from the dynamical sampling theory.  We study the following

\bigskip

{\bf Main Problem.}
\\
Let $\Lambda$ be a uniformly discrete subset of $\R^n$ and let  ${G_{\alpha}(x)}$ be a collection of functions parametrized by $\alpha \in I$.
What assumptions should be imposed on the spatial set $\Lambda$, index set $I$, and functions $G_{\alpha}$ to enable the recovery of every band-limited signal $f$ from its space-time samples $\{f \ast G_{\alpha}(\lambda)\}_{\lambda \in \La, \alpha \in I}$?

\bigskip

For signals $f$ from a Paley-Wiener space $PW_{\sigma}$ (see the definition below) it means that the inequalities
\begin{equation}\label{sspw}
    D_1 \|f\|^2_{2} \leq  \sum \limits_{\lambda \in \La} \int \limits_{I} |f * G_{\alpha}(\lambda)|^2 \, d\alpha \leq D_2 \|f\|^2_{2} \quad \text{for every } f \in PW_{\sigma}
\end{equation}
are true with some constants $D_1$ and $D_2$. Here, as usual, $\|\cdot\|$ denotes the $L^2$-norm.

Recall that a set $\Lambda = \{\lambda_k\} \subset\R^n$ is called uniformly discrete\footnote{Sometimes, the term uniformly separated is used.} (u.d.) if  
$$\delta(\Lambda) := \inf\limits_{\substack{\lambda \neq \lambda' \\ \lambda, \lambda' \in \La}}{|\lambda - \lambda'|}  > 0.$$

The constant $\delta(\Lambda)$ is called the separation constant of $\Lambda$.

In the one-dimensional setting, this problem appears in particular in connection with tasks of mathematical physics. Several examples are presented in \cite{g}. One of them is the initial value problem for the heat equation
\begin{equation}\label{ivp_eq}
\frac{\partial}{\partial \alpha} u(x, \alpha) =
\sigma^2 \frac{\partial^2 u}{\partial x^2} (x,\alpha), \quad \quad \sigma \neq 0,\,  x \in \R, \, \alpha > 0,
\end{equation}
with initial condition
\begin{equation}\label{init_ivp_eq}
u(x,0) = f(x).
\end{equation}
It is well-known that the solution is given by the formula
\begin{equation}
    u(x,\alpha) = f \ast g_{\alpha} (x)= \int\limits_{\R^n} g_{\alpha} (x - y) f(y) dy,
\end{equation}
where  $ g_{\alpha}(x) = \frac{1}{\sqrt{(4 \pi \alpha \sigma)}} \exp \left( -\frac{x^2}{4 \alpha \sigma} \right).$ Note that Main Problem applied to equation~(\ref{ivp_eq}) provides the reconstruction of initial function $f$ from the states $\{u(\lambda, \alpha)\}_{\lambda \in \La, \alpha \in I}$.

A variant of Main problem for the one-dimensional setting was considered by Aldroubi et al. in \cite{g}.
In particular, it was established that unlike the classical sampling setting, the assumptions that should be imposed on the set $\Lambda$ to solve the Main Problem cannot be expressed in terms of some density of $\Lambda$, see Example~4.1 in \cite{g}. More precisely, one may construct a set with an arbitrarily small density that provides stable reconstruction of the initial signal.
Also in that paper, it was shown that for the solution of Main Problem we have to require $\Lambda$ to be relatively dense.

In the one-dimensional setting, for a large collection of kernels, a solution of Main Problem was presented in  \cite{UZ}: It turns out the stable recovery from the samples on $\Lambda$ is possible if and only if $\Lambda$ is not (in a certain sense) “close” to an arithmetic progression.

It seems natural to extend the results of \cite{g} and \cite{UZ} to the multi-dimensional situation. Below we focus on the two-dimensional variant of the problem for the case of Gaussian kernel

$$G_{\alpha}(x) = e^{- \alpha_1 x_1^2 - \alpha_2 x_2^2}, \quad \quad \alpha = (\alpha_1, \alpha_2) \in I, \, I= (a,b) \times (c,d) \subset \R^2_+, \quad  x = (x_1, x_2) \in \R^2. 
$$
Our approach is similar to the one in \cite{UZ}. However,  this problem is considerably more involved than the one in the one-dimensional setting. One needs to apply some additional ideas.
See Section~\ref{section_mult_and_op} for some remarks on cases dimension higher than $2$.

We pass to the description of the geometry of the sets $\Lambda$ that solve the planar Main problem.

\begin{Def}
A curvilinear lattice in $\R^2$ defined by three vectors
$$t=(t_1,t_2) \in \R^2, \quad  \xi=(\xi_1,\xi_2)\in\R^2, \quad \text { and }  \quad r=(r_1,r_2) \in \R^2,\, r_1^2 + r_2^2 = 1,$$ is the set of all vectors $\lambda=(\lambda_1,\lambda_2)\in\R^2$ satisfying $$
l_{t,\xi,r}:=\{\lambda \in \R^2 \; \Big| \; r_1\cos(\lambda_1\xi_1+\lambda_2\xi_2+t_1)=r_2\cos(-\lambda_1\xi_1+\lambda_2\xi_2+t_2)\}.
$$
\end{Def}
{
\begin{figure}[h!]\label{counter_ex}
\center{\includegraphics[height=7cm]{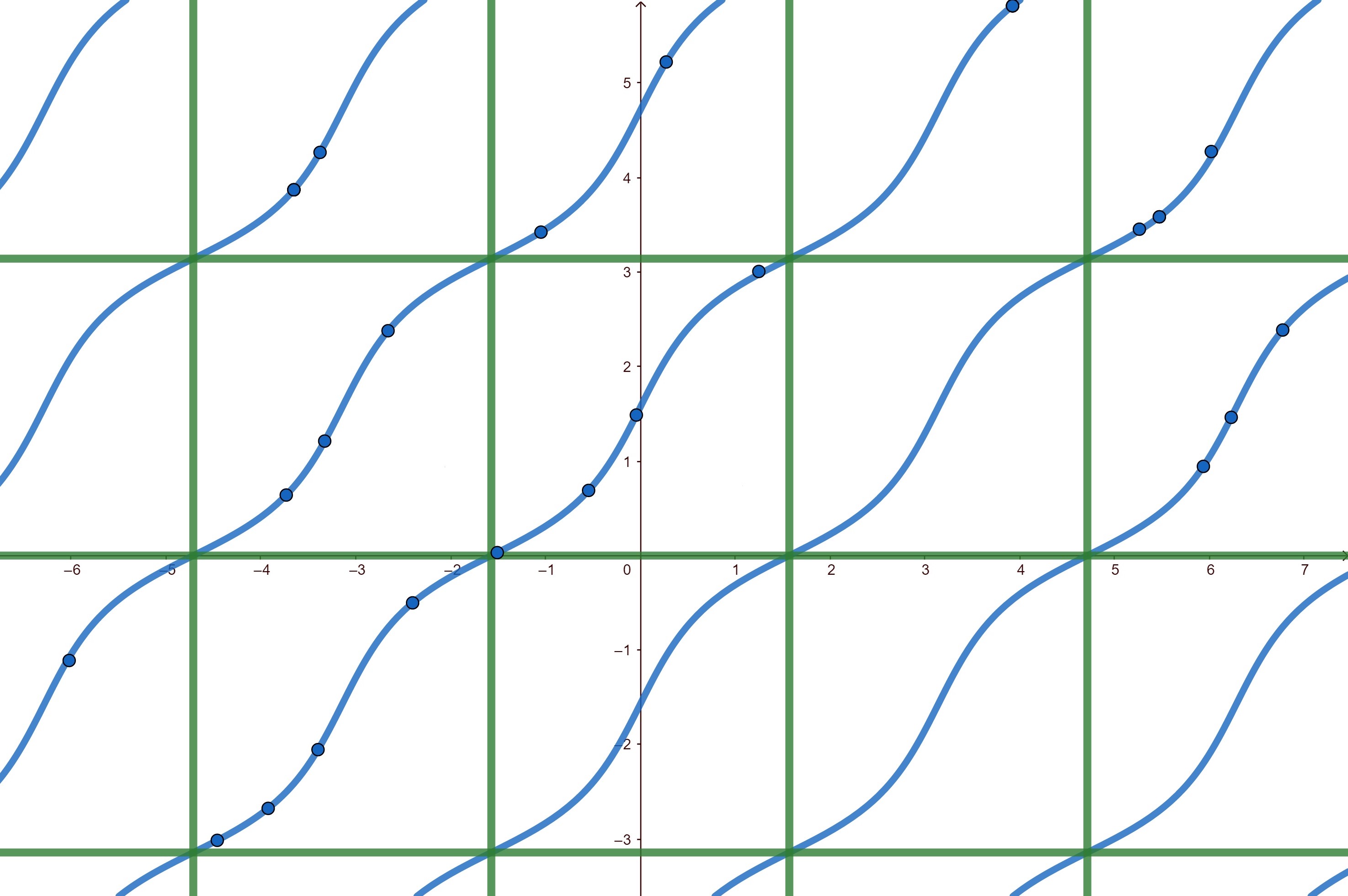}}
\caption{Curvilinear lattice defined by $\cos (x+y) = 3 \cos(y-x)$}
\label{img_curve_lat}
\end{figure}
}
The blue curves on Figure~\ref{img_curve_lat} correspond to the curvilinear lattice $l_{t, \xi, r}$ with $t=(0,0), \xi = (1,1)$, and  $r = (1/\sqrt{10}, 3/\sqrt{10})$. 

In what follows the notation $W(\Lambda)$ stands for the collection of all weak limits of translates of a uniformly discrete set $\Lambda$, see the definition in Section~\ref{section_notations}.

\bigskip\noindent{\bf Condition (A)}: A uniformly discrete set $\La=\{\lambda=(\lambda_1,\lambda_2)\}\subset\R^2$ satisfies condition (A) if  every set $\La^\ast\in W(\La)$ is not empty and does not lie on any lattice
$l_{t,\xi,r}$.

\begin{remark}
A Delone set is a set that is both uniformly discrete and relatively dense.
In particular, it is easy to check that every set that satisfies condition (A) is a Delone set.
\end{remark}

We denote by $PW^2_{\sigma}$ the space of square integrable on $\R^2$ functions with spectrum supported in the square $[-\sigma, \sigma]^2$, i.e. 
$$
PW^2_{\sigma} = \{f \in L^2(\R^2) \,\, \big| \,\, \supp \hat{f} \subset [-\sigma, \sigma]^2 \},
$$
where
$$
\hat{f}(\xi_1, \xi_2) = \int\limits_{\R^2} e^{-i (\xi_1  x_1 + \xi_2 x_2)} f(x_1, x_2) \,  dx_1 dx_2.
$$

Now, we are ready to formulate the main result.

\begin{theorem}\label{Main_theorem}
Given a u.d. set $\La\subset\R^2$ and a rectangle $I = (a,b)\times(c, d)$ with $0 < a < b < \infty$, $0 < c < d < \infty$.
The following statements are equivalent:

{\rm (i)} For every $\sigma>0$ there are positive constants $D_1 = D_1(\sigma, I, \Lambda)$ and $D_2 = D_2(\sigma, I, \Lambda)$ such that {\rm(\ref{sspw})} holds true.

{\rm (ii)} $\La$ satisfies condition {\rm (A)}.
\end{theorem}

The paper is organized as follows. In Section~\ref{section_notations} we give all necessary definitions and fix some notations.
As it was mentioned above, we employ the approach from \cite{UZ} and divide the solution into two parts. We start with solving Main Problem for the Bernstein spaces $B_{\sigma}$ and prove an analogue of Theorem~\ref{Main_theorem} in Section~\ref{sec_bernstein}. In Section~\ref{pw_section} we investigate the connection between the sampling with Gaussian kernel in the Paley-Wiener and Bernstein spaces. We also prove the main result in Section~\ref{pw_section}. The remarks on multi-dimensional cases and some open problems that puzzle us are placed in Section~\ref{section_mult_and_op}.

\section{Notations and preliminaries}\label{section_notations}

In the present paper, we deal with signals that belong to the Bernstein and Paley-Wiener spaces.
Since we investigate Main Problem simultaneously for all bandwidth parameters, we may consider only the functions with the spectrum supported in squares. This leads us to
\begin{Def}
Given a positive number $\sigma$, we denote by $B_{\sigma}$ the space of all entire functions $f$ in $\C^2$ satisfying the estimate
\begin{equation}\label{exp_est}
|f(z)|\leq C e^{\sigma(|y_1|+|y_2|)},\quad z=(z_1,z_2)\in \C^2, \quad z_j=x_j+i y_j\in \C, \, j=1,2,
\end{equation}
where the constant $C=C(f)$ depends only on $f$.
\end{Def}
It is well-known that $B_{\sigma}$ consists of the bounded continuous functions that are the inverse Fourier transforms of tempered distributions supported on the square $[-\sigma, \sigma]^2$.  We refer the reader to \cite{levin} for more information about Bernstein spaces.

For $1\le p < \infty$ we may define the Paley-Wiener spaces by the formula
$$PW^p_{\sigma} = B_{\sigma} \cap L^p(\R^2) $$
or equivalently
$$
PW^p_{\sigma} = \{f \in L^p(\R^2) \,\, \big| \,\, \supp \hat{f} \subset [-\sigma, \sigma]^2 \}.
$$

Following \cite{B} (see also Chapter 3.4 in \cite{ou1}, \cite{JNR}, and \cite{RUZ}), we introduce auxiliary

\begin{Def}
Let $\{\Lambda_k\}$ and $\Lambda$ be u.d. subsets of $\R^n$, satisfying $\delta(\Lambda_k) \ge \delta > 0, \, k \in \N$. We say that the sequence 
$\{\Lambda_k\}$ converges weakly to $\Lambda$ if for every large $ R > 0$ and small $\varepsilon > 0$ there exists such $N = N(R, \varepsilon)$ that 
$$
\Lambda_{k} \cap (-R, R)^n \subset \Lambda + (-\varepsilon, \varepsilon)^n,
$$
$$
\Lambda \cap (-R, R)^n \subset \Lambda_{k} + (-\varepsilon, \varepsilon)^n.
$$
for all $k \ge N$.
\end{Def}

\begin{Def}
By $W(\Lambda)$ we denote all weak limits of the translates $\Lambda_k:=\Lambda - x_k$, where $\{x_k\} \subset \R^n$ is an  arbitrarily bounded or unbounded sequence.
\end{Def}

We supply these definitions with several examples concerning the condition (A).
\begin{ex}
To construct the set that does not satisfy condition (A), one may consider the following perturbation of the rectangle lattice:
$$
\Lambda = \left\{ \left(2 \pi n + \frac{1}{2^{m^2+n^2}}, 2 \pi m + \frac{1}{2^{|m| + |n|}} \right), \quad m,n \in \Z \right\}.
$$
Taking any sequence $\{x_k\} \subset \R^2$ such that $|x_k| \to \infty$, one may check by the definition that the sequence $\Lambda - x_k$ weakly converges to the set $\Lambda = \left\{ \left(2 \pi n, 2 \pi m  \right), \quad  m,n \in \Z \right\}$, which, clearly, lies in $l_{t, \xi, r}$ with $t = (0,0)$, $\xi = (1,1)$, and $r = \left( \frac{1}{\sqrt{2}}, \frac{1}{\sqrt{2}} \right)$.
\end{ex}


The following example is inspired by the papers \cite{JNR} and \cite{RUZ}, which were dedicated to the solution of planar mobile sampling problems.
\begin{ex}
Set 
$$D_{\Z} = \left\{ (x,y) \subset \R^2 \quad | \quad x^2 + y^2 = 4 \pi^2 k^2, \,\, k \in \Z  \right\},$$
i.e. $D_{\Z}$ is a collection of concentric equidistant circles with center $(0,0)$.
Now, we may take as $\Lambda$ the u.d. set located on the circles $D_{\Z}$.
One may check (see the proofs in {\rm \cite{JNR}} and {\rm\cite{RUZ}}) that the weak limits of translates for every unbounded sequence $\{x_n\}$ for $D_{\Z}$ lies on the parallel lines. The argument is based on the simple observation that the traces of translated circles in the rectangle $[-R,R]^2$ (for a fixed $R>0$) are getting closer and closer to the lines as the value $|x_n|$ increases. Moreover, the distance between these lines is $2 \pi k$. For instance, one may take $x_n = (0,2\pi n)$ and pass to a weak limit $\Lambda - x_n \to \Lambda^{\ast}$ to obtain that $\Lambda^{\ast} \subset l_{t,\xi,r}$ with $t = (0,0)$, $\xi = (1,1)$, and $r = \left( \frac{1}{\sqrt{2}}, \frac{1}{\sqrt{2}} \right)$.
\end{ex}

Below we will use the simple fact that for every sequence $x_k$ there is a subsequence $x_{k_j}$ such that $\Lambda - x_{k_j}$ converges weakly.



Throughout this paper we will adopt the following notations: 
\begin{itemize}
\item  Let $x \in \R^n, \, y \in \R^n, \, n \in \N$. Define $|x| := \sqrt{x_1^2 + \dots + x_n^2}$. Notation $x \cdot y$ stands for the scalar product of vectors $x$ and $y$.  
    
\item  Set $B_r(x) := \{y \in \R^n \, : \, |x-y| < r\},\,\,$ where $x \in \R^n$ and $r>0$. 

\item Given $\lambda=(\lambda_1, \dots, \lambda_n)$ and $f \in L^{\infty}(\R^n)$, we set
$$
 f_\lambda(x):=f(x-\lambda) = f(x_1 -\lambda_1, \dots, x_n -  \lambda_n).
$$

\item By $|A|$ we denote the $n$-th dimensional Lebesgue measure of a set $A \subset \R^n$.    
    
\item By $C$ we denote different positive constants.

\end{itemize}

Basically, we will focus on two-dimensional case. It is convenient to fix the following notations.

\begin{itemize}

\item Given a point $x=(x_1,x_2)\in\R^2$, denote $\tilde{x}:=(-x_1,x_2)$.

\item A symmetrization operator $S$ is defined by the formula
$$
Sf(x):=f(x)+f(\tilde{x})+f(-\tilde{x})+f(-x), \quad f \in L^{\infty}(\R^2).
$$

\item Set $\T := \{|x| = 1, \, x \in \R^2\}$.

    

\end{itemize}

\section{Sampling with Gaussian kernel in Bernstein spaces}\label{sec_bernstein}

An analogue of Theorem~\ref{Main_theorem} for the Bernstein spaces is as follows:

\begin{theorem}\label{ST_B_theorem}
Given a u.d. set $\La\subset\R^2$ and $I = (a,b)\times(c,d)$ with $0 < a < b < \infty, 0 < c < d < \infty$.
The following statements are equivalent:

{\rm (i)} For every $\sigma>0$ there is a constant $K=K(\sigma)$ such that 
$$
\|f\|_\infty\leq K\sup_{\alpha\in I} \sup_{\lambda \in \La} \|f\ast G_\alpha\|_\infty \quad \text{ for every } f\in B_\sigma.
$$

{\rm (ii)} $\La$ satisfies condition {\rm (A)}.
\end{theorem}

Above, as usual, $\|\cdot\|_{{\infty}}$ denotes the sup-norm 
$$
\|f\|_{\infty} := \sup_{x \in \R^2} |f(x)|.
$$


\subsection{Proof of Theorem \ref{ST_B_theorem}, Part I}

(ii) $\Rightarrow$ (i).  In what follows we assume that (i) is not true. We have to show that $(ii)$ fails, i.e. there is a set $\La^\ast \in W(\La)$ such that it lies on some curvilinear lattice. 
The proof is divided into $5$ steps. For the convenience of the reader, we will briefly describe them here and then pass 
to the argument.

In Step 1, using the standard Beurling technique, we find  $\Lambda^{\ast} \in W(\Lambda)$ and $g \in B_{\sigma}$ such that $g \ast G_{\alpha}$ vanishes on $\Lambda^{\ast}$ for every $\alpha \in I$. Our next step is to show that $Sg_{\lambda} = 0$ for every $\lambda \in \La^{\ast}$. In Step 3 we prove that $\La^{\ast}$ lies on some curvilinear lattice under the assumption that $g \in L^2(\R^2)$. In Steps 4 and 5, using some approximation technique, we show how to get rid of the requirement that $g$ is square integrable.

1.   Due to the assumption made one can find a sequence of Bernstein functions $f_n\in B_\sigma$ satisfying
$$
\|f_n\|_\infty=1,\quad \|f_n\ast G_\alpha|_\La\|_\infty\leq 1/n.
$$We may then introduce a sequence of functions
$$g_n(z):=f_n(z-x(n)), \quad z=(z_1,z_2), \quad x(n)=(x_1(n),x_2(n)),$$
where $x(n)$ are chosen so that $|f_n(x(n))| > 1 - \frac{1}{n}, \, n \in \N$. Then we have $$\|g_n\|_{\infty}=1 \quad \text{and} \quad \left\|g\ast G_{\alpha}|_{\Lambda+x(n)}\right\|_{\infty} \leq 1/n,\, n\in \N.$$ Using the compactness property of Bernstein space (see, e.g., \cite{ou1}, Proposition~2.19), we may assume that sequence $g_n$ converges (uniformly on compacts in $\C^2$) to some function $g\in B_\sigma$.
Moreover, passing if necessary to a subsequence, we may assume that the translates $\Lambda+x(n)$ converge weakly to some u.d. set $\Lambda^\ast$. Of course, we may assume that $\Lambda^{\ast}$ is non-empty. Otherwise, we have arrived at contradiction with condition (A).
Clearly, $g$ satisfies
\begin{equation}\label{g}
\|g\|_\infty=1,  \text{ and for every } \alpha\in I: \quad  g\ast G_\alpha|_{\La^\ast}=0, \quad  \La^\ast\in W(\La).
\end{equation}

For a point $z = (z_1, z_2) \in \C^2$ we consider its complex conjugate point $\bar{z} = (\bar{z_1},\bar{z_2})$.

Consider the decomposition $g(z) = \varphi(z) + i \psi(z)$, where
$$
\varphi(z):=\frac{g(z)+\overline{g(\bar{z})}}{2}, \quad 
\psi(z):=\frac{g(z)-\overline{g(\bar{z})}}{2i}.
$$
Then  $\varphi$ and $\psi$ are real (on $\R^2$) entire functions satisfying (\ref{exp_est}).
Thereby, functions $\varphi$ and $\psi$ belong to $\B_\sigma$, and since the kernel $G_{\alpha}$ takes only real values on $\R^2$, we have 
$(\varphi\ast G_{\alpha})(\lambda)=0$ and $(\psi\ast G_{\alpha})(\lambda)=0$ for every $\lambda \in \La^{\ast}$. Thus, we can continue the argument assuming that $g$ is a real-valued function.

2. Recall that the notations $\tilde{x}$ and $Sf$ were introduced in Section~\ref{section_notations}.

\begin{lemma}\label{l1}
Assume a function $g\in B_\sigma$ satisfies {\rm(\ref{g})}. Then for every $\lambda\in \La^\ast$ the equality
\begin{equation}\label{sg}
 S g_\lambda(x)=0
\end{equation}
holds for a.e. $ x \in \R^2$.
\end{lemma}
\noindent{\bf Proof}. 
Without loss of generality, we may assume that $\lambda = (0,0)$ and $I=\left(\frac{1}{2},1\right)^2$.
Observe that 
\begin{equation}\label{SgG}
(S g \ast G_{\alpha})(0,0) = 4 (g * G_{\alpha})(0,0) = 0 \quad \text{ for every } \alpha \in I.
\end{equation}

Set
$$
h(x_1, x_2) := Sg(x_1, x_2) \exp\left\{ - \frac{x_1^2 + x_2^2}{4} \right\} \quad \text{ and } \quad I_+ := \left(\frac{1}{2}, \frac{3}{4}\right)^2.
$$
Clearly, $h \in L^2(\R^2)$ and it is even in variables $x_1$ and $x_2$.
Moreover, using~(\ref{SgG}), one can check that $(h \ast G_{\alpha})(0,0) = 0$ for any $\alpha \in I_+$.

For every multi-index $m = (m_1, m_2) \in \mathbb{N}^2$ and $u \in I_+$ we have
$$
\frac{\partial^{m}}{\partial u_1^{m_1} \partial u_2^{m_2}} \int\limits_{\R^2} h(x_1,x_2) \exp\left\{ - u_1 x_1^2 - u_2 x_2^2 \right\} dx_1 dx_2 = 0.
$$
In particular, $h$ is orthogonal to every monomial $x_1^{\alpha_1} x_2^{\alpha_2}$ with even indexes $\alpha_1$ and $\alpha_2$ in the weighted space $L^2\left(\R^2, \exp\left\{ - \frac{1}{2}  (x_1^2+x_2^2) \right\} \right)$. 
Moreover, since $h$ is even in any variable, from the symmetry, we see that $h$ is orthogonal to every polynomial in this space. To finish the proof we use the completeness property of multi-dimensional analogues of Hermite polynomials. More precisely, we invoke Theorem~3.2.18 from \cite{ort_pol} to deduce $h = 0$. Consequently, $Sg(x) =0$ for every $x \in \R^2$, and the lemma follows.
$\Box$

3. We will need a simple technical
\begin{lemma}\label{l2}
Given a function $F\in L^2(\R^2)$ such that its inverse Fourier transform $f$ is a real function. Then
$$
S f_\lambda (x)=2\int\limits_{\R^2}\cos (x \cdot t)\, \mbox{\rm Re}\left(e^{i\lambda\cdot t}F(t)+e^{i\tilde{\lambda}\cdot t}F(\tilde{t})\right)dt.
$$
\end{lemma}

\noindent{\bf Proof}. Indeed, we may write
$$
f(x)=\mbox{\rm Re}\int\limits_{\R^2}e^{i x\cdot t}F(t)dt.
$$Therefore,
$$
Sf_\lambda (x)=\mbox{\rm Re}\int\limits_{\R^2}\left(e^{i x\cdot t}+e^{i\tilde{x}\cdot t}+e^{-i x\cdot t}+e^{-i\tilde{x}\cdot t}\right)e^{i\lambda\cdot t}F(t)dt=$$$$2\mbox{\rm Re}\int\limits_{\R^2}(\cos (x\cdot t)+\cos(\tilde{x}\cdot t))e^{i\lambda\cdot t}F(t)dt=2\mbox{\rm Re}\int\limits_{\R^2}\cos(x\cdot t)\left(e^{i\lambda\cdot t}F(t)+e^{i\tilde{\lambda}\cdot t}F(\tilde{t})\right)dt ,
$$ which proves the lemma. $\Box$

\medskip

3. If we additionally assume that $g\in L^2(\R^2)$, the result follows from the next statement.
\begin{lemma}\label{l3}
Assume  $g\in L^2(\R^2)$. Then $\La^\ast$ lies on some curvilinear lattice.
\end{lemma}

\noindent{\bf Proof}. Denote by $G$ the inverse Fourier transform of $g$. Recall that $g$ is real, whence
$$
G(t) = \overline{G(-t)} \quad \text{ and } \quad G(\tilde{t}) = \overline{G(-\tilde{t})}.
$$
Denote by $U(t) = \mbox{\rm Re}\left(e^{i\lambda\cdot t}G(t)+e^{i\tilde{\lambda}\cdot t}G(\tilde{t})\right)$.
Since
$$
2U(t) =  2\mbox{\rm Re}\left(e^{i\lambda\cdot t}G(t)+e^{i\tilde{\lambda}\cdot t}G(\tilde{t})\right) = \left(e^{i\lambda\cdot t}G(t)+e^{i\tilde{\lambda}\cdot t}G(\tilde{t})\right) + \overline{\left(e^{i\lambda\cdot t}G(t)+e^{i\tilde{\lambda}\cdot t}G(\tilde{t})\right)} =  ,
$$
$$
e^{i\lambda\cdot t}G(t)+e^{i\tilde{\lambda}\cdot t}G(\tilde{t}) + 
e^{-i\lambda\cdot t}G(-t)+e^{-i\tilde{\lambda}\cdot t}G(-\tilde{t}),
$$
we deduce that $U(t) = U(-t)$.
Combining this observation with Lemmas \ref{l1} and \ref{l2}, we see that equality
$$
 \mbox{\rm Re}\left(e^{i\lambda\cdot t}G(t)+e^{i\tilde{\lambda}\cdot t}G(\tilde{t})\right)=0
$$
holds for a.e. $t\in\R^2$ and for every $\lambda \in \La^{\ast}$.

Recall that $G=0$ a.e. outside $(-\sigma,\sigma)^2$. For every $\epsilon>0$, find a real Schwartz function $F_\epsilon$ whose support lies on $[-\sigma,\sigma]^2$ satisfying $\|G-F_\epsilon\|_2<\epsilon$. Then
\begin{equation}\label{FGeps}
\|F_\epsilon\|_2\geq \|G\|_2-\epsilon
\end{equation}
and
\begin{equation}\label{ref_est}
\left(\,\, \int\limits_{\R^2} \left|\mbox{\rm Re}\left(e^{i\lambda\cdot t}F_\epsilon(t)+e^{i\tilde{\lambda}\cdot t}F_\epsilon(\tilde{t})\right)\right|^2 d t \right)^{1/2}<2\epsilon,\quad \lambda\in\La^\ast.
\end{equation}
Using these inequalities, one can check that there are a point $t_\epsilon$, $t_{\epsilon} \in [-\sigma, \sigma]^2$ and constants $C$ and $c$ depending only on $\sigma$ 
such that
\begin{equation}\label{t_eps_ineq}
|F_\epsilon(t_\epsilon)|> C \quad \text{ and } \quad  \Bigl|\mbox{\rm Re}\left(e^{i\lambda\cdot t_\epsilon}F_\epsilon(t_\epsilon)+e^{i\tilde{\lambda}\cdot t_\epsilon}F_\epsilon(\tilde{t_\epsilon})\right)\Bigr|< c \epsilon |F_\epsilon(t_\epsilon)|,    
\end{equation}
for every small enough $\epsilon$.
Indeed, by Plancherel theorem, we have $\|G\|_2 = \|g\|_2$. Since $\|g\|_{\infty} = 1$ and $g \in B_{\sigma}$, using Bernstein inequality, we deduce $\|G\|_2 = \|g\|_2 \ge C(\sigma)$. Now, assuming that for all $t$,  the inequalities \eqref{t_eps_ineq} do not hold true, by integration with respect to variable $t$ and using the estimate \eqref{FGeps}, for sufficiently small $\epsilon$ we arrive at
$$
\int\limits_{\R^2} \left|\mbox{\rm Re}\left(e^{i\lambda\cdot t}F_\epsilon(t)+e^{i\tilde{\lambda}\cdot t}F_\epsilon(\tilde{t})\right)\right|^2 d t \ge c^2 \epsilon^2 \|F_{\epsilon}\|^2_2 \ge \frac{c^2}{2} \epsilon^2 \|G\|^2_2 \ge \frac{c^2}{2} \epsilon^2 C^2(\sigma),
$$
which contradicts to estimate \eqref{ref_est} when $c > 2 \sqrt{2}/C(\sigma)$.

Write
$$
F_\epsilon(t_\epsilon)=:R_\epsilon e^{iu_\epsilon},\quad F_\epsilon(\ti_\epsilon)=:r_\epsilon e^{iv_\epsilon}.
$$Then we get
$$
|R_\epsilon|>C,\quad \left|R_\epsilon\cos(\lambda\cdot t_{\epsilon}+u_\epsilon)+r_\epsilon\cos(\tilde{\lambda}\cdot t_\epsilon+v_\epsilon)\right|< c \epsilon R_\epsilon, \ \lambda\in\La^\ast.
$$

Then normalizing we arrive at
$$
\left|\frac{R_\epsilon}{\sqrt{R_\epsilon^2+r_\epsilon^2}}\cos(\lambda\cdot t_{\epsilon}+u_\epsilon)+\frac{r_\epsilon}{\sqrt{R_\epsilon^2+r_\epsilon^2}}\cos(\tilde{\lambda}\cdot t_\epsilon+v_\epsilon)\right|< c\epsilon.$$
Clearly, we may assume that $u_{\epsilon} \in [0, 2 \pi]$ and $v_{\varepsilon} \in [0, 2 \pi]$. Recall that $t_{\epsilon} \in [-\sigma, \sigma]^2$ and, of course,
$$
\frac{R_\epsilon}{\sqrt{R_\epsilon^2+r_\epsilon^2}} \in [0,1] \quad \text{and} \quad \frac{r_\epsilon}{\sqrt{R_\epsilon^2+r_\epsilon^2}} \in [0,1].
$$
Taking $\epsilon = \frac{1}{n}$ and passing if necessary to a subsequence, we deduce that $\Lambda$ lies on some curvilinear lattice.

4. In what follows we assume that \begin{equation}\label{gg} g\in B_{\sigma} \setminus L^2(\R^2).\end{equation}

For $\epsilon > 0$ we set
$$
h_{\epsilon}(\xi):=\frac{\sin(\epsilon \xi)}{ \epsilon \xi}, \quad \Phi_{\epsilon}(x_1, x_2) := h_{\epsilon}(x_1) h_{\epsilon}(x_2),  \quad
\text{and} \quad
\delta_\epsilon:=\|g \Phi_{\epsilon}\|_2^{-1/2}.
$$

The next statement easily follows from (\ref{gg}). 
\begin{lemma}\label{l}
We have $\delta_\epsilon\to0$ as $\epsilon\to 0.$
\end{lemma}
We skip the simple proof.

Let us introduce auxiliary functions$$
\varphi_\epsilon(x):=\delta_\epsilon \Phi_{\epsilon}(x),\quad g_{\epsilon}(x):=g(x)\varphi_\epsilon(x), \quad x\in \R^2.
$$By Lemma \ref{l},
\begin{equation}\label{ge}
\|g_{\epsilon}\|_2=1/\delta_\epsilon\to \infty,\quad \epsilon\to 0.
\end{equation}


\begin{lemma}
For every $\lambda\in\La^\ast$ satisfying $|\lambda|< 1/\sqrt{\delta_\epsilon}$ we have
$$
\|S (g_\epsilon)_\lambda \|_2 \leq C\sqrt\delta_\epsilon.
$$
\end{lemma}

\noindent{\bf Proof}. 
By Lemma~\ref{l1}, $S g_{\lambda} = 0$. Since the function $\varphi_\epsilon$ is even with respect to each variable, we have
$$S g_{\lambda} \varphi_{\epsilon}(x) = (Sg(\cdot-\lambda)\varphi_\epsilon(\cdot))(x) =
$$
$$g(x - \lambda) \varphi_{\epsilon}(x) + g(\tilde{x} - \lambda) \varphi_{\epsilon}(\tilde{x}) + g(-x - \lambda) \varphi_{\epsilon}(-x) + g(-\tilde{x} - \lambda) \varphi_{\epsilon}(-\tilde{x}) = 0.$$
Hence,
$$
\left|S (g_\epsilon)_\lambda (x)\right| = \left|S (g \varphi_\epsilon)(\cdot-\lambda) (x)\right| = \left|S (g \varphi_\epsilon)(x-\lambda) - Sg(\cdot-\lambda)\varphi_\epsilon)(x) \right| \leq
$$
$$
\left|g(x - \lambda)(\varphi_\epsilon(x-\lambda)-\varphi_\epsilon(x))\right|+
\left|g(\tilde{x}-\lambda)(\varphi_\epsilon(\tilde{x}-\lambda)-\varphi_\epsilon(\tilde{x}))\right| + 
\left|g(-x - \lambda)(\varphi_\epsilon(-x-\lambda)-\varphi_\epsilon(-x))\right|+$$
$$+\left|g(-\tilde{x}-\lambda)(\varphi_\epsilon(-\tilde{x}-\lambda)-\varphi_\epsilon(-\tilde{x}))\right|.$$

Below we focus on the estimate of the first term at the right hand-side of the inequality above. The remaining terms admit the same estimate.

Write $\lambda=(\lambda_1,\lambda_2)$. Observe that
$$
 \left|\varphi_\epsilon(x-\lambda)-\varphi_\epsilon(x)\right| \leq \delta_\epsilon\Bigl(\left|h_{\epsilon}(x_1 - \lambda_1) - h_{\epsilon}(x_1)\right| \left|h_{\epsilon}(x_2 - \lambda_2)\right|+\left|h_{\epsilon}(x_2-\lambda_2)-h_{\epsilon}(x_2)\right| |h_{\epsilon}(x_1)|\Bigr).
 $$

For $j = 1,2$ using the Cauchy-Schwartz inequality, we have
 $$\left(\,\,\int\limits_{\R} |h_{\epsilon}(x_j-\lambda_j)-h_{\epsilon}(x_j)|^2 dx_j \right)^{1/2} =
 \left(\,\,\int\limits_{\R}
 \left| \int\limits_{0}^{\lambda_j}h_{\epsilon}'(x_j-u) du \right|^2 dx_j \right)^{1/2} \leq C |\lambda_j| \|h_{\varepsilon}'\|_{2}.$$
%
One may check that $\|h_{\epsilon}\|_2=C/\sqrt{\epsilon}$ and $\|h'_{\epsilon}\|_2=C\sqrt{\epsilon}$.
Since $\|g\|_\infty=1$ and $|\lambda| \le 1/\sqrt{\delta_{\epsilon}}$, we arrive at
$$
\|S (g_\epsilon)_\lambda \|_2 \leq C \left(\,\,\int\limits_{\R^2} \left|\varphi_\epsilon(x-\lambda)-\varphi_\epsilon(x)\right|^2 dx\right)^{1/2} \leq C \delta_{\epsilon} |\lambda| \|h_{\epsilon}\|_{2} \|h'_{\epsilon}\|_{2}  \leq C\sqrt\delta_\epsilon.
$$
That finishes the proof. 
$\Box$

5. Denote by $G_\epsilon:=\widehat{g\varphi_\epsilon}$. Then $G_\epsilon\in L^2(\R)$  vanishes a.e. outside some square 
$(-\sigma^\ast,\sigma^\ast)^2$ (it is easy to check that one may take $\sigma^\ast=\sigma+\epsilon)$.

Using Lemma 5, for $|\lambda| \leq 1/ \sqrt{\delta_{\epsilon}}$ we get
$$
\left(\,\, \int\limits_{\R^2}\left|\mbox{\rm Re}\left(e^{i\lambda\cdot x}G_\epsilon(x)+e^{-i\lambda\cdot\tilde{x}}G_\epsilon(\tilde{x})\right)\right|^2 dx\right)^{1/2} \leq C\sqrt\delta_\epsilon.
$$
On the other hand, by (\ref{ge}), $\|G_\epsilon\|_2\geq C$, for all small enough $\epsilon$. 

To finish the proof, we proceed as in the proof of Lemma 3.

\subsection{Proof of Theorem \ref{ST_B_theorem}, Part II}

(i) $\Rightarrow$ (ii). We will argue by contradiction. 
Assume that for every $\sigma > 0$ there is a constant $K=K(\sigma)$ such that
$$
\|f\|_\infty\leq K\sup_{\alpha\in I} \sup_{\lambda \in \La} \|f\ast G_\alpha\|_\infty,\quad f\in B_\sigma,
$$ but condition (ii) is not satisfied, i.e. there exists some $\Lambda' \in W(\La)$ such that $\Lambda'$ lies on some curvilinear lattice. 
Clearly, to come to the contradiction it suffices to construct for every $\varepsilon > 0$ a  function $f = f_{\varepsilon}$ such that
\begin{equation}\label{ST2_f_prop}
    \|f\|_{\infty} \ge C, \quad \sup\limits_{\alpha \in I} \sup\limits_{\lambda \in \La} |f * G_{\alpha}(\lambda)| \le C\varepsilon,
\end{equation}
and $f \in \B_{\sigma^{\ast}}$ for some fixed $\sigma^{\ast}$.
 
Again, let us provide a brief description of the proof.
We divide the proof into $4$ steps. First, we build a function $g$ such that $g \ast G_{\alpha}$ vanishes on $\Lambda'$ for every $\alpha \in I$. A slight modification of $g$ provides a function $f$, which satisfies~(\ref{ST2_f_prop}). To verify the second estimate in ~(\ref{ST2_f_prop}) we split the set $\Lambda$ into the sets $ \Lambda_I = \Lambda \cap P$ and $\Lambda_O = \Lambda \cap (\R^2 \setminus P)$ for an appropriate rectangle $P$. In the steps 3 and 4, we show that $f$ satisfies the relations~(\ref{ST2_f_prop}) for $\lambda \in \La_O$ and $\lambda \in \Lambda_I$ respectively. 

Now we pass to the proof.

1. By our assumption, there exist $\Lambda' \in W(\Lambda)$, $\xi \in \R^2$, $(t_1,t_2) \in \R^2$,  and $(r_1, r_2) \in \T$ such that for every $\lambda' \in \La'$ the equality
\begin{equation}\label{curve_th_eq}
     r_1\cos(\lambda' \cdot \xi - t_1) - r_2  \cos(\tilde{\lambda'} \cdot \xi - t_2) = 0 
\end{equation}
holds. 
Set
$$
   g(x) =  r_1 \cos( \xi \cdot x + t_1) - r_2  \cos(  \tilde{\xi} \cdot x  + t_2).   
$$
Clearly, $g \in B_{\sigma}$ for $\sigma = |\xi|$. Next, we will show that symmetrization of the function $g_{\lambda'}$ vanishes for every $\lambda' \in \ \La'$.
\begin{lemma}
    The equality
$$
    S g_{\lambda'}(x) = 0 
$$
holds for every $x \in \R^2$ and $\lambda' \in \La'$.
\end{lemma}
\noindent{\bf Proof}.
After some simple calculations, we have
$$
S g_{\lambda'} (x) = r_1 \, \re \left(e^{i ( t_1 - \xi \cdot \lambda' )} S(e^{i \xi \cdot x}) \right) - r_2 \, \re \left(e^{i (t_2 -  \tilde{\xi} \cdot \lambda')} S(e^{i \tilde{\xi} \cdot x}) \right),$$
where we, as usual, apply symmetrization operator $S$ with respect to variable $x$. Clearly, $S(e^{i \tilde{\xi} \cdot x}) = S(e^{i \xi \cdot x}) = 2 (\cos (\xi \cdot x) + \cos ( \tilde{\xi} \cdot x))$. Thus, using ~(\ref{curve_th_eq}), we have
$$
S g_{\lambda'} (x) = 2 (\cos (\xi \cdot x) + \cos ( \tilde{\xi} \cdot x)) \left(r_1 \, \re e^{i (t_1 - \lambda' \cdot \xi)} + r_2 \, \re e^{i (t_2 - \tilde{\lambda'} \cdot \xi)} \right) = 0.
$$
$\Box$

Consequently, for every $\lambda' \in \Lambda'$ and $\alpha \in I$ we have
\begin{equation}
     g \ast G_{\alpha}(\lambda') = 0,
\end{equation}
since $G_{\alpha}$ is even in every variable.

2. Fix small $\varepsilon > 0$ and take large $R=R(\varepsilon) > 0$ (we will specify its value later). 
Recall that $\La' \in W(\Lambda)$. In particular, that means that one can find $v=(v_1, v_2) = v(R, \varepsilon) \in \R^2$ such that inside the square $[-R,R]^2$, the set  $\Lambda - v$ is "close"\, to $\Lambda'$:
\begin{equation}\label{lambda_neighbour}
    \text{for every } \lambda \in  \La \cap (v + (-R,R)^2) \text{ there is } \lambda' \in \La' \cap(-R,R)^2: \quad \dist(\lambda - v, \lambda') \le \varepsilon.
\end{equation}
Set $P = [v_1-R,v_1+R]\times[v_2-R,v_2+R]$ and consider the decomposition
$$
    \Lambda = \Lambda_I \cup \Lambda_O := \left(\Lambda\cap P \right) \cup \left(\Lambda\cap (\R^2 \setminus P) \right).
$$
Consider
$$
    \Phi_{\varepsilon}(t) = \Phi_{\varepsilon}(t_1, t_2) = \frac{\sin(\varepsilon t_1)}{\varepsilon t_1} \frac{\sin(\varepsilon t_2)}{\varepsilon t_2}.
$$
We define the function $f$ by the formula
$$
    f(x) = \Phi_{\varepsilon}(x-v) g(x-v), \quad x\in \R^2, \, v \in \R^2.
$$
Clearly, $\|f\|_{\infty}  \geq C$, and it suffices to show that $|f \ast G_{\alpha}(\lambda)| \leq C \varepsilon$ for every $\lambda \in \Lambda$. 
We will estimate the value $|f \ast G_{\alpha}(\lambda)|$ for $\lambda \in \Lambda_I$ and $\lambda \in \Lambda_O$ separately. 

3. Assume that $\lambda \in \Lambda_O$. We may choose $R = R(\varepsilon) = \frac{1}{\varepsilon^2}$.
Set $U = U_1 \times U_2 =  [-\sqrt R, \sqrt{ R}]^2.$
For $s\in U$ we have
\begin{equation}\label{sq_est}
    |f(\lambda - s)| \le \frac{\|g\|_{\infty}}{\varepsilon^2 |\lambda_1 - s_1 - v_1| |\lambda_2 - s_2-v_2|}  \le \frac{\|g\|_{\infty}}{\varepsilon^2 |R-\sqrt{R}|^2} \le C  \varepsilon^2 \|g\|_{\infty}, 
\end{equation}
since $|\lambda_1 - v_1| \ge R$ and $|\lambda_2-v_2| \ge R$.
Next, it is easy to check that
\begin{equation}\label{strip_est}
      J:=\int\limits_{\R}\int\limits_{\R \setminus U_1} G_{\alpha}(s_1,s_2) \, ds_1 ds_2+ \int\limits_{\R} \int\limits_{\R \setminus U_2} G_{\alpha}(s_1,s_2) \, ds_2 ds_1 \le C \varepsilon. 
\end{equation}

Now, to estimate $f \ast G_{\alpha} (\lambda)$ for $\lambda \in \La_O$, we write 
$$
    |f\ast G_{\alpha}(\lambda)| \le \int\limits_{U} |f(\lambda -s)| G_{\alpha} (s) \, ds + \int\limits_{\R} \int\limits_{\R \setminus U_1} |f(\lambda -s)| G_{\alpha} (s) \, ds + \int\limits_{\R} \int\limits_{\R \setminus U_2} |f(\lambda -s)| G_{\alpha} (s) \, ds.
$$
Applying $\|f\|_{\infty} \le 1$ and estimates~(\ref{sq_est}) and~(\ref{strip_est}), we arrive at

$$
    |f\ast G_{\alpha}(\lambda)| \le C\varepsilon^2 \int\limits_{U} G_{\alpha} (s) \, ds +  J \le C \varepsilon.
$$

4. Now, assume that $\lambda \in \Lambda_I$. Take $\lambda' \in \Lambda'$, satisfying condition~(\ref{lambda_neighbour}) corresponding to $\lambda$, i.e. $\dist(\lambda - v, \lambda' ) < \varepsilon$.
Since $g\ast G_{\alpha}(\lambda') = 0$, we may write
$$
    f\ast G_{\alpha}(\lambda) = \int\limits_{\R^2} f(\lambda-s) G_{\alpha}(s) ds + \Phi_{\varepsilon}(\lambda-v) \int\limits_{\R^2} g(\lambda'-s) G_{\alpha}(s) ds =  
$$
$$
\int\limits_{\R^2} \big( \Phi_{\varepsilon}(\lambda-s-v)\left(g(\lambda-s-v) - g(\lambda'-s) \right) +  g(\lambda'-s)\left( \Phi_{\varepsilon}(\lambda-v-s) - \Phi_{\varepsilon}(\lambda-v) \right) \big)  G_{\alpha}(s)  ds. 
$$
Set
$$
    H_1:=\left|\Phi_{\varepsilon}(\lambda - s - v) -  \Phi_{\varepsilon}(\lambda - v) \right|,
$$
$$
    H_2:=\left|g(\lambda - s - v) - g(\lambda' - s)\right|.
$$
Clearly,
\begin{equation}\label{inside_rect}
    |f\ast G_{\alpha}(\lambda)| \le \int\limits_{\R^2}
    \left(H_1 |g(\lambda' - s)| +
    H_2 \left| h_{\varepsilon}(\lambda - v) \right|  \right)
    G_{\alpha}(s) ds. 
\end{equation}
By Bernstein inequality and relation (\ref{lambda_neighbour}), we have
\begin{equation}\label{H_estimate}
    |H_1| \le \varepsilon (|s_1| + |s_2|), \quad |H_2| \le \varepsilon \|g'\|_{\infty} .
\end{equation}
Combining estimates (\ref{inside_rect}) and (\ref{H_estimate}) together, we obtain
$$
    |f\ast G_{\alpha}(\lambda)| \le C \varepsilon \int\limits_{\R^2}
    (\|g'\|_{\infty} + |s_1|+|s_2|)
    G_{\alpha}(s_1, s_2) ds_1 ds_2 \le C \varepsilon 
$$
that finishes the proof. $\Box$

The following statement easily follows from Theorem~\ref{ST_B_theorem}.
\begin{lemma}\label{p_lemma}
Assume $\Lambda$ and $I$ satisfy the assumptions of Theorem~\ref{ST_B_theorem} and condition $(i)$ is fulfilled. 
Then for every $\sigma > 0$ there is a constant $C$ such that
\begin{equation}\label{infty_p_estimate}
    \|f\|_{\infty}^2 \le C \int\limits_{I} \sup \limits_{\lambda \in \La} | f \ast G_{\alpha} (\lambda)|^2 \, d\alpha
\end{equation}
for every $f \in PW^2_{\sigma}$.
\end{lemma}

\section{Sampling with Gaussian kernel in Paley-Wiener spaces}\label{pw_section}
\subsection{Auxiliary statements}\label{aux_b_pw}
Recall that our aim is to describe the geometry of sets $\La \subset \R^2$ that for every $f \in PW^2_{\sigma}$ the estimates 
\begin{equation}\label{pw_estimates}
 D_1 \|f\|^2_{2} \leq  \sum \limits_{\lambda \in \La} \int \limits_{I} |f * G_{\alpha}(\lambda)|^2 \, d\alpha \leq D_2 \|f\|^2_{2}, 
\end{equation}
hold with some constants $D_1$ and $D_2$ independent on $f$.

\subsubsection{Bessel-type inequality}

We start with showing that the right hand-side of (\ref{pw_estimates}) follows easily from classical sampling results for u.d. set $\Lambda$.

\begin{proposition}\label{revers_proposition}
Assume $\Lambda$ is a u.d. set,
$I = (a,b) \times (c,d)$, where $0<a<b<\infty, 0<c<d<\infty$.
Then there is a constant $D_2 = D_2(I,\Lambda)$ such that
$$
     \sum \limits_{\lambda \in \La} \int \limits_{I} |f * G_{\alpha}(\lambda)|^2 \, d\alpha  \le D_2 \|f\|^2_{2}, 
$$ for every $f \in PW^2_{\sigma}$.
\end{proposition}
\noindent{\bf Proof}.
Recall the {\it Bessel inequality} for Paley-Wiener spaces:  
if $\Lambda$ is a u.d. subset of $\R^d$ then there is a constant $M=M(\La, \sigma)$ such that  
 \begin{equation}\label{PP_est}
     \sum\limits_{\lambda \in \La}|g(\lambda)|^2 \le M \|g\|^2_{2}
 \end{equation}
for every $g \in PW^2_{\sigma}$, see \cite{Y}, Chapter 2, Theorem 17.

Note that convolution with Gaussian Kernel $G_{\alpha}$ keeps the function in Paley-Wiener space. Using Young's convolution inequality and Bessel inequality, one can find a constant $D_2$ such that for every $f \in PW^2_{\sigma}$ the estimate 
$$
    \sum \limits_{\lambda \in \La} \int \limits_{J} |f * G_{\alpha}(\lambda)|^2 \, d\alpha \le C |J| \|f \ast G_{\alpha}\|^2_{2} \le C \|f\|_{2}^2\|G_{\alpha}\|_{1}^2 \le D_2 \|f\|_{2}^2
$$
is true. 

That finishes the proof of proposition. $\Box$

\subsubsection{Auxiliary functions}
In what follows we need some auxiliary functions with special properties. These functions should belong to Paley-Wiener spaces, have a large $L^2$-norm with a small $L^2$-norm of the gradient. 
Now, we specify these requirements. 

\medskip\noindent{\bf Condition (B)}: Let $\varepsilon$ be a small positive parameter. A family of functions $\{\Phi_{\varepsilon}\}$ satisfies condition ($B$) if 
{
\begin{enumerate}
    \item[($\beta_1$)] $\Phi_{\varepsilon}(0, 0) = 1, \quad \|\Phi_{\varepsilon}\|_{\infty} = 1;$ 
    \item[($\beta_2$)] $\Phi_{\varepsilon} \in PW^2_{\varepsilon}$;
    \item[($\beta_3$)] $\|\Phi_{\varepsilon}\|_{2} \to \infty$ as $\varepsilon \to 0$;
    \item[($\beta_4$)] $\| \nabla \Phi_{\varepsilon}\|_{2} \to 0$ as $\varepsilon \to 0$.
\end{enumerate}
}
Next, we provide a few examples to illustrate some additional difficulties that occur in the multi-dimensional setting. Then we present an example of functions $\Phi_{\varepsilon}$ that satisfy condition (B).

\begin{ex}\label{sinc_ex}
Let us return to the one-dimensional case. Consider
$$
\Phi_{\varepsilon} (x) = \frac{\sin (\varepsilon x)}{\varepsilon x}.
$$
Observe that functions $\Phi_{\varepsilon}$ satisfy an analogue of condition (B) in the  one-dimensional setting.
Clearly, $\Phi_{\varepsilon} (0) = 1, \|\Phi_{\varepsilon}\|_{\infty} = 1$, and $\Phi_{\varepsilon} \in PW^2_{\varepsilon}$.
One may easily check that $$\| \Phi_{\varepsilon}\|_{2} \le C \varepsilon^{-1/2} \quad \text{and} \quad \| \Phi'_{\varepsilon}\|_{2} \le C \varepsilon^{1/2}.$$
These relations prove the one-dimensional analogues of $(\beta_3)$ and $(\beta_4)$.
\end{ex}

The passage from Bernstein to Paley-Wiener spaces and back in \cite{UZ} was based on the properties of the functions in Example~\ref{sinc_ex}. One may try to construct functions $\Phi_{\varepsilon}$ that satisfy condition (B) in the two-dimensional setting in the following natural way.

\begin{ex}\label{ex_phi_2_sinc}
Consider the function $\Phi_{\varepsilon}$ defined by the formula
$$
\Phi_{\varepsilon} (x, y) = \frac{\sin (\varepsilon x)}{\varepsilon x} \frac{\sin (\varepsilon y)}{\varepsilon y}.
$$

It is clear that conditions $(\beta_1)$ and $(\beta_2)$ are true. Property $(\beta_3)$  follows from
$$
\|\Phi_{\varepsilon}\|_{2}\le C \varepsilon^{-1}.
$$
However, one may easily check that $\|\nabla \Phi_{\varepsilon}\|_{2}$ does not converge to zero as $\varepsilon \to 0$. 

\end{ex}
However, in two-dimensional setting it is still possible to construct functions that satisfy condition (B).
Now, we pass to the construction. 

\begin{lemma}\label{DS_lemma}
Assume $\varepsilon > 0$.  
There exist functions $\Psi_{\varepsilon}$ such that
\begin{enumerate}
    \item[$(P1)$] $\supp \Psi_{\varepsilon} \subset B_{\varepsilon}(0), \quad \Psi_{\varepsilon} \ge 0,$
    \item[$(P2)$] $C_1 \le \int\limits_{\R^2} \Psi_{\varepsilon}(x) \, dx  \le C_2, \quad  0 < C_1 \le C_2 < \infty,$
    \item[$(P3)$] $\|\Psi_{\varepsilon}\|_{2} \ge  \frac{C}{\varepsilon^{3/4}},$
    \item[$(P4)$] $\left(\int\limits_{\R^2} |\Psi_{\varepsilon}(x)|^2 |x|^2 \, dx \right)^{1/2} \le \frac{C}{\sqrt{ \log \frac{1}{\varepsilon}}}.$
\end{enumerate}
\end{lemma}
\noindent{\bf Proof}. Fix small $0 < \varepsilon < 1$ and denote the integer part of $\log \frac{1}{\varepsilon}$ by $m$.
For integers $n$ from $[m, 2m]$ we set  $a_n = 2^{2n}/n$.
Next, we define the function $\Psi_{\varepsilon}$ layer by layer by the formula
$$
\Psi_{\varepsilon}(x) = a_n, \quad x \in B_{2^{-n}}(0) \setminus B_{2^{-n-1}}(0).
$$
For $|x|>\varepsilon$ and  $|x| < \frac{\varepsilon^2}{2}$ we set $\Phi_\varepsilon(x)=0$.
Note that the area of the ring $B_{2^{-n}}(0) \setminus B_{2^{-n-1}}(0)$ is equal to $\frac{3 \pi}{4} 2^{-2n}$.

Clearly, $\Psi_{\varepsilon}$ satisfy (P1). To verify (P2) we write
$$
\int\limits_{\R^2} \Psi_{\varepsilon}(x) \, dx = \frac{3\pi}{4} \sum\limits_{n = m}^{2m}  2^{-2n}  a_n = C
\sum\limits_{n = m}^{2m} \frac{1}{n}.
$$
Note that the right-hand side of this equation can be estimated with some fixed positive constants from above and below by 
$$ \int\limits_{\log \frac{1}{\varepsilon}}^{2\log \frac{1}{\varepsilon}} \frac{1}{t} \, dt = \log 2.$$  Thus, condition (P2) follows. Next, we have
$$
\|\Psi_{\varepsilon}\|^2_{2} = \frac{3 \pi}{4} \sum\limits_{n = m}^{2m}  2^{-2n} a^2_n \ge C \sum\limits_{n = m}^{2m} \frac{2^{2n}}{n^2} \ge C \int\limits_{\log \frac{1}{\varepsilon}}^{2\log \frac{1}{\varepsilon}} 2^{2t} t^{-2}\, dt \ge \frac{C}{\varepsilon^2 \log^2 \frac{1}{\varepsilon}} \ge \frac{C}{\varepsilon^{3/2}}, 
$$
and $(P3)$ follows. The estimate
$$
\int\limits_{\R^2} |\Psi_{\varepsilon}(x)|^2 |x|^2 \, dx \le C  \sum\limits_{n = m}^{2m} 2^{-4n} a^2_n \le C  \sum\limits_{n = m}^{2m} \frac{1}{n^2} \le \frac{C}{ \log \frac{1}{\varepsilon}}
$$
implies $(P4)$ that finishes the proof. $\Box$

\begin{corollary}\label{B_22_holds}
There exist functions $\Phi_\varepsilon$ satisfying condition {\rm(B).}
\end{corollary}
\noindent{\bf Proof}. Denote by $c_{\Psi} = \int\limits_{\R^2} \Psi_{\varepsilon}(x) \, dx$. By $(P_2)$, $c_{\Psi}$ is positive, finite, and separated from zero.
Now, we may define $\Phi_\varepsilon$ as the Fourier transform of $\Psi_{\varepsilon}$ with a proper normalization:
$$
\Phi_{\varepsilon}(x) = \frac{1}{c_{\Psi}} \int_{\R^2}e^{-i x\cdot t}\Psi_{\varepsilon}(t) \,dt.
$$
The property $(\beta_2)$ follows from $(P_1)$. Due to $\Psi_{\varepsilon} \ge 0$ and normalization condition $(\beta_1)$ is fulfilled. Relations $(P_3)$ and $(P_4)$ imply estimates $(\beta_3)$ and $(\beta_4)$ respectively. $\Box$


\subsection{From Bernstein to Paley-Wiener spaces and back}
To prove Theorem~\ref{Main_theorem}, we will use the following statement, which describes the connection between sampling in Paley-Wiener and Bernstein spaces.

\begin{theorem}\label{PW_B_p_theorem}
Let $\La$ be a u.d. set in $\R^2$, $I = (a,b) \times (c,d), 0 < a < b < \infty, 0 < c < d < \infty,$ and $\sigma'>\sigma>0$. 

{\rm (i)} Assume the inequality
\begin{equation}\label{th_bs_ss}
\|f\|_\infty\leq K\sup_{\alpha\in I} \sup_{\lambda \in \La} \|f\ast G_\alpha\|_\infty \quad \text{  for all  } f\in B_{\sigma'}
\end{equation}
holds with some constant $K = K(\sigma', \Lambda)$. Then there exists a constant $D_1 = D_1(\sigma, \Lambda)$ such that
\begin{equation}\label{th_pw_ss}
D_1 \|f\|^2_{2} \leq  \sum \limits_{\lambda \in \La} \int \limits_{I} |f * G_{\alpha}(\lambda)|^2 \, d\alpha \quad \text{  for every  }  f \in PW^2_{\sigma}
\end{equation}
is true.

{\rm (ii)} Assume that {\rm (\ref{th_pw_ss})} holds with some constant $D_1 = D_1(\sigma', \Lambda)$  for all $f \in PW_{\sigma'}$. Then there is a constant $K = K(\sigma', \Lambda)$ such that  {\rm (\ref{th_bs_ss})} is true for every $f \in B_{\sigma}$. 
\end{theorem}
\begin{remark}
For a similar result for space sampling see {\rm \cite{OUM}.}
\end{remark}
\begin{remark}
In this theorem we do not need to require $I$ to be a rectangle. One may take $I=(a,b) \subset \R$ with $0< a < b < \infty$. In such a case by $G_{\alpha}(x)$ we mean $G_\alpha(x_1,x_2) = e^{-\alpha (x_1^2+x_2^2)}$ and $d \alpha$ is a standard one-dimensional Lebesgue measure.
\end{remark}

The proof of Theorem~\ref{PW_B_p_theorem} is similar to the proof of Theorem~3  in the paper~\cite{UZ}.
We provide the argument for statement $(i)$ and leave the proof of $(ii)$ to the reader.
The functions $\Phi_{\varepsilon}$ that satisfy condition (B) play a crucial role in our argument.

\medskip

\noindent{{\bf Proof of Theorem~\ref{PW_B_p_theorem}. }}
Take $\varepsilon > 0$ such that $\sigma + \varepsilon < \sigma'$.
By our assumption, for every $q \in B_{\sigma}$ the estimate 
\begin{equation}
    \|q\|_{\infty} \le C \sup\limits_{\alpha \in I} \sup\limits_{ \lambda \in \La} |q * G_{\alpha}(\lambda)|, 
\end{equation}
is true and our aim is to prove~(\ref{th_pw_ss}).



Consider functions $\Phi_{\varepsilon}$ satisfying condition (B). 

Using $(\beta_1)$, we get
\begin{equation}\label{eq_PWB1}
    \|f\|^2_{2} = \int\limits_{\R^2} |f(x)|^2 dx \le \int\limits_{\R^2} \sup\limits_{t \in  \R^2} \left| \Phi_{\varepsilon}(x - t)  f(t) \right|^2 dx.
\end{equation}
Note that $q(t) := \Phi_{\varepsilon}(x - t) f(t) \in \B_{\sigma + \varepsilon}$, and we can apply Lemma \ref{p_lemma} to obtain
\begin{equation}\label{eq_PWB2}
    |q(t)|^2 \le C \int\limits_I \sup\limits_{\lambda \in \Lambda} \left |\,\int\limits_{\R^2} G_{\alpha}(\lambda - s) \Phi_{\varepsilon}(x - s) f(s) \, ds \right|^2 d \alpha,
\end{equation}
where the constant $C$ does not depend on $t$.
To provide the estimate from above we may replace $\sup\limits_{\lambda \in \La}$ by $\sum\limits_{\lambda \in \La}$, and switch the order of integration and summation:
\begin{equation}\label{eq_PWBp}
    \|f\|^2_{2} \le C
    \sum\limits_{\lambda \in \Lambda} \int\limits_I  \int\limits_{\R^2}  \left |\, \int\limits_{\R^2} G_{\alpha}(\lambda - s) \Phi_{\varepsilon}(x - s) f(s) \, ds \right|^2  \, dx \, d \alpha.
\end{equation}
Denote by
\begin{equation}\label{eqPWBJ1}
    Y_1 = \left| \Phi_{\varepsilon}(x - \lambda) \,\int\limits_{\R^2} G_{\alpha}(\lambda - s)  f(s) \, ds \right|^2,
\end{equation}
\begin{equation}\label{eqPWBJ2}
    Y_2 = \left|\,\int\limits_{\R^2} G_{\alpha}(\lambda - s)  \left(   \Phi_{\varepsilon} (x - \lambda)  -  \Phi_{\varepsilon}(x - s) \right) f(s) \, ds \right|^2.
\end{equation}
Using the inequality $|a + b|^2 \le C(|a|^2 + |b|^2)$, we deduce from (\ref{eq_PWBp}), (\ref{eqPWBJ1}), and (\ref{eqPWBJ2}) that
\begin{equation}\label{J_1+J_2}
    \|f\|^2_{2} \le  C \sum\limits_{\lambda \in \Lambda} \int\limits_I \int\limits_{\R^2}    \left(Y_1 + Y_2\right)  \, dx \, d\alpha. 
\end{equation}
Next, we estimate the terms with $Y_1$ and $Y_2$ separately.
The value of $\sum\limits_{\lambda \in \Lambda} \int\limits_I \int\limits_{\R^2}Y_1\, dx \, d\alpha $ is majorized by  
\begin{equation}\label{J_1_est}
    \sum\limits_{\lambda \in \Lambda} \int\limits_I  \left(\,\, \int\limits_{\R^2}
    | \Phi_{\varepsilon}(x - \lambda)|^2 \, dx \right) |(f \ast G_{\alpha})(\lambda)|^2  \, d\alpha  \le  
    \|\Phi_{\varepsilon}\|_{2}^2
    \int\limits_I \sum\limits_{\lambda \in \Lambda} |(f \ast G_{\alpha})(\lambda)|^2 d \alpha.
\end{equation}
The inequalities for the second term are more complicated. Set
$$
    H(x; \lambda,s) = |\Phi_{\varepsilon}(x - \lambda) -  \Phi_{\varepsilon}(x - s)|.
$$
We start with the observation 
$$
      H(x; \lambda,s)  \le  \left| \int\limits_{s_1}^{\lambda_1} \frac{\partial \Phi_{\varepsilon}}{\partial x}(x - u_1, y - \lambda_2) du_1 \right| + \left| \int\limits_{s_2}^{\lambda_2} \frac{\partial \Phi_{\varepsilon}}{\partial y}(x - s_1, y - u_2) du_2 \right|.   
$$
Using Cauchy-Schwarz inequality, we write 
$$
     H^2(x; \lambda,s)  \le C \Bigg((\lambda_1 - s_1)  \int\limits_{s_1}^{\lambda_1} \left| \frac{\partial \Phi_{\varepsilon}}{\partial x}(x - u_1, y - \lambda_2)\right|^2 du_1 +
     $$
     $$
     (\lambda_2 - s_2)  \int\limits_{s_2}^{\lambda_2} \left| \frac{\partial \Phi_{\varepsilon}}{\partial x}(x - s_1, y - u_2)\right|^2 du_2  \Bigg).
$$
Thereby, for $\lambda = (\lambda_1, \lambda_2)$ and $s = (s_1, s_2)$ we get
$$
    \int\limits_{\R^2} H^2(x; \lambda,s) \, dx \le  C |s_1 - \lambda_1| |s_2 - \lambda_2| \|\nabla \Phi_{\varepsilon}\|_{2}^2,
$$
whence
\begin{equation}\label{H2_est}
    \|H(\,\cdot\,; \lambda, s)\|^2_{2} \le  C |s- \lambda|^{2}  \|\nabla \Phi_{\varepsilon}\|_{2}^2.
\end{equation}
Now, we return to the estimation of the term with $Y_2$ in the formula~(\ref{J_1+J_2}).
Applying Cauchy--Schwarz inequality, we arrive at
$$
\sum\limits_{\lambda \in \Lambda} \, \int\limits_{\R^2} Y_2 \, dx  =  
\sum\limits_{\lambda \in \Lambda}\,\,  \int\limits_{\R^2} \left|\,\, \int\limits_{\R^2} f(s) G_{\alpha}(\lambda - s) H(x; \lambda,s) \, ds \right|^2 dx \le 
$$
$$
\sum\limits_{\lambda \in \Lambda}  \,\int\limits_{\R^2} \left( \,\, \int\limits_{\R^2} |f(s)|^2 G_{\alpha}(\lambda - s) ds \, \int\limits_{\R^2} G_{\alpha}(\lambda - s) H^2(x; \lambda,s) ds \right) dx. 
$$
With estimate (\ref{H2_est}) in hand, we continue 
$$
\sum\limits_{\lambda \in \Lambda}\,  \int\limits_{\R^2} Y_2 \, dx \le \sum\limits_{\lambda \in \Lambda}  \left( \,\, \int\limits_{\R^2} |f(s)|^2 G_{\alpha}(\lambda - s) ds  \,\int\limits_{\R^2} G_{\alpha}(\lambda - s) \|H^2(\cdot; \lambda,s)\|_{2}  ds  \right) \le
$$
$$
\le C  \|\nabla \Phi_{\varepsilon}\|_{2}^2 \sum\limits_{\lambda \in \Lambda}  \left( \,\, \int\limits_{\R^2} |f(s)|^2 G_{\alpha}(\lambda - s) ds  \,\int\limits_{\R^2} G_{\alpha}(\lambda - s) |s - \lambda|^2  ds  \right).
$$
Clearly,
\begin{equation}\label{JJJ}
      \int\limits_{\R^2} G_{\alpha}(\lambda - s) |s - \lambda|^2  ds  \le C,
\end{equation}
and since $\Lambda$ is a u.d. set, we have
\begin{equation}\label{G_sum_est}
    \sum \limits_{\lambda \in \La}
    G_{\alpha}(\lambda - s) \le C.
\end{equation}
Using relations~(\ref{JJJ}) and (\ref{G_sum_est}), we finish the estimate of the term with $Y_2$:
\begin{equation}\label{J_2_est}
    \sum\limits_{\lambda \in \Lambda}\; \int\limits_{I}  \int\limits_{\R^2} Y_2\, dx \, d \alpha \le  C |I| \|\nabla \Phi_{\varepsilon}\|_{2}^2  \|f\|_{2}^2. 
\end{equation}
Combining (\ref{J_1+J_2}), (\ref{J_1_est}), and (\ref{J_2_est}) together, we get
\begin{equation}
    \|f\|^2_{2} \le C_1 |I| \|\nabla \Phi_{\varepsilon}\|_{2}^2 \|f\|_{2}^2 +  C_2 \|\Phi_{\varepsilon}\|_{2}^2 
    \int\limits_I \sum\limits_{\lambda \in \Lambda} |(f \ast G_{\alpha})(\lambda)|^2 d \alpha.
\end{equation}

To finish the proof we invoke properties $(\beta_3)$ and $(\beta_4)$. Indeed, taking sufficiently small $\varepsilon > 0$  we make the first summand less than $\frac{\|f\|^2_{2}}{2}$, and (\ref{th_pw_ss}) follows. $\Box$

\subsection{Proof of Theorem~\ref{Main_theorem}}
Now, we are ready to prove the main result.

(i)  $\Rightarrow$ (ii) \medskip
Assume that for the set $\Lambda$ condition (i) is satisfied. In particular,  for any $\sigma > 0$ inequality~(\ref{th_pw_ss}) is true for every $f\in PW^2_{\sigma}$ with constant $D_1$ depending on $\sigma$.
Then, Theorem~\ref{PW_B_p_theorem} implies that for every $\sigma > 0$ inequality~(\ref{th_bs_ss}) holds true for every $f \in B_{\sigma}$ with constant $K$ depending on $\sigma$. By Theorem~\ref{ST_B_theorem}, we deduce that $\Lambda$ satisfy condition (A).

(ii)  $\Rightarrow$ (i)
Assume that condition (ii) is fulfilled.
Recall that Proposition~\ref{revers_proposition} ensures that the right hand-side estimate in (\ref{sspw}) holds with some universal constant. Thus, it suffices to verify that inequality~(\ref{th_pw_ss}) is true for every $\sigma > 0$ and every $f \in PW^2_{\sigma}$ with some constant $D_1 = D_1(\sigma)$. By our assumption and Theorem~\ref{ST_B_theorem}, the inequality~(\ref{th_bs_ss}) is true for every $\sigma > 0$ and every $f \in B_{\sigma}$ with a constant $K$ depending only on $\sigma$. Applying Theorem~\ref{PW_B_p_theorem}, we see that~(\ref{th_pw_ss}) holds true for every $\sigma > 0$  and $f \in PW^2_{\sigma}$ with a constant $D_1$ depending only on $\sigma$. Thus, condition (i) is true. That finishes the proof. $\Box$


\section{Remarks}\label{section_mult_and_op}

First, we would like to note that Theorems~\ref{Main_theorem},  \ref{ST_B_theorem}, and \ref{PW_B_p_theorem} remain true for a collection of kernels that satisfy some additional assumptions similar to conditions $(\beta) - (\theta)$ in \cite{UZ}.


Second, one may check that our approach provides a solution to the Main Problem for the Bernstein spaces $B_{[-\sigma, \sigma]^n}$ for the Gaussian kernels in $\R^n$, with an index set $I = \prod\limits_{i=1}^n[a_i,b_i]$, and  $\Lambda$ a u.d. subset of $\R^n$. One may therefore formulate an analogue of Theorem~\ref{ST_B_theorem} in multi-dimensional setting.
However, the passage to Paley-Wiener spaces faces obstacles similar to those discussed
in Section~\ref{pw_section}, Example~\ref{ex_phi_2_sinc}.
We note that the frame inequalities for some continuous frame $\{e_x\}_{x \in X}$
\begin{equation}\label{frame_eq}
    D_1 \|f\|^p_{p} \leq  \int \limits_{X} | \langle f, e_{x} \rangle|^p dx \, \leq D_2 \|f\|^p_{p}
\end{equation}
(compare with \eqref{sspw} and consider $X = I \times \Lambda $ equipped with the measure $dx$ --- product of $n$-dimensional Lebesgue measure on $I$ and counting measure on $\Lambda$)
typically hold true for all range of Banach spaces $(X, \|\cdot\|_p), \, 1 \le p < \infty$ simultaneously if the frame ${e_x}$ has a sufficiently good localization, see \cite{A0}, \cite{Gr_fr}, and \cite{Gr_fr2}. However, in our setting we didn't manage to prove the analog of main result for all $p \in [1, \infty)$ in multi-dimensional case $d>2$.

On the other hand, using our approach, one may check that for every $n > 2$ there are a number $p(n)$ and functions $\Phi_{\varepsilon}$ such that for $p \ge p(n)$ we have 
 $$\|\Phi_{\varepsilon}\|_{p} \to \infty,  \quad \| \nabla \Phi_{\varepsilon}\|_{p} \to 0 \quad \text{ as } \varepsilon \to 0.$$
Thus, a modification of the proof of Theorem~\ref{PW_B_p_theorem} leads to the solution of the Main problem for $PW^p_{[-\sigma, \sigma]^n}$ spaces with $p \ge p(n)$. 




\section{Acknowledgements}
I am grateful to A.~Ulanovskii for stimulating discussions and to D.~Stolyarov for the proof of Lemma~\ref{DS_lemma}. 

\noindent Ilya Zlotnikov\\
St. Petersburg Department of V.A. Steklov \newline
	    Mathematical Institute of the Russian Academy of Sciences \newline
	    27 Fontanka, St. Petersburg 191023, Russia,\\
	    zlotnikk@rambler.ru

\end{document}